%% file: 2003-17.tex
\documentclass{gtart}
\input gtoutput

\lognumber{300}
\volumenumber{7}\papernumber{17}
\volumeyear{2003}
\pagenumbers{615}{639}
\received{16 January 2003}
\revised{17 October 2003}
\accepted{21 September 2003}
\published{22 October 2003}
\proposed{Robion Kirby}
\seconded{Tomasz Mrowka, Cameron Gordon}

\usepackage{amsmath,amsfonts,epsf,amscd}

\newcommand{\Tors}{\mathrm{Tors}}
\newcommand{\HF}{HF}

\newtheorem{theorem}{Theorem}[section]
\newtheorem{prop}[theorem]{Proposition}
\newtheorem{cor}[theorem]{Corollary}

\newtheorem{lemma}[theorem]{Lemma}

\newcommand{\Q}{\mathbb{Q}}
\newcommand{\R}{\mathbb{R}}

\newcommand{\C}{\mathbb{C}}

\newcommand{\Z}{\mathbb{Z}}

\newcommand{\Zmod}[1]{\Z/{#1}\Z}

\newcommand{\Image}{\mathrm{Im}}

\newcommand{\cm}{\cdot}

\newcommand{\ModSWfour}{\mathcal{M}}
\newcommand{\ModFlow}{\ModSWfour}

\newcommand{\SpinC}{{\mathrm{Spin}}^c}

\newcommand\Wedge{\Lambda}

\newcommand\Hom{\mathrm{Hom}}

\newcommand\abuts\Rightarrow
\newcommand\Sym{\mathrm{Sym}}

\newcommand\mCP{{\overline{\mathbb{CP}}}^2}

\newcommand\relspinc{\underline{\spinc}}

\newcommand\Filt{\mathcal F}

\newcommand\x{\mathbf x}

\newcommand\z{\mathbf z}

\newcommand\y{\mathbf y}

\newcommand\ModSphere{\ModFlow\left({\mathbb S}\longrightarrow 
\Sym^{g-1}(\Sigma_{1})\times \Sym^2(\Sigma_{2})\right)}
\newcommand\ModSpheres\ModSphere
\newcommand\CF{CF}

\newcommand\CFa{\widehat{CF}}
\newcommand\CFp{\CFb}
\newcommand\CFm{\CF^-}

\newcommand\HFp{\HFb}

\newcommand\HFm{\HF^-}
\newcommand\CFinf{CF^\infty}
\newcommand\HFinf{HF^\infty}
\newcommand\CFb{CF^+}
\newcommand\HFa{\widehat{HF}}
\newcommand\HFb{HF^+}

\newcommand\UnparModSp{\widehat \ModSp}
\newcommand\UnparModFlow\UnparModSp
\newcommand\Mod\ModSp

\newcommand{\spinc}{\mathfrak s}

\newcommand\Real{\mathrm Re}

\newcommand\ModMaps{\mathcal M}
\newcommand\ModSp\ModMaps

\newcommand\Ta{{\mathbb T}_{\alpha}}
\newcommand\Tb{{\mathbb T}_{\beta}}

\newcommand\alphas{\mbox{\boldmath$\alpha$}}

\newcommand\betas{\mbox{\boldmath$\beta$}}

\newcommand\Dual{\mathcal D}
\newcommand\Duality\Dual

\newcommand\Seif{\Sigma}

\newcommand\Char{\mathrm{Char}}
\newcommand\Conc{\mathrm{Conc}}

\newcommand\spincrel\relspinc

\newcommand\CFK{CFK}
\newcommand\HFK{HFK}

\newcommand\CFKa{\widehat\CFK}

\newcommand\CFKinf{\CFK^{\infty}}

\newcommand\HFKa{\widehat\HFK}

\newcommand\BasePt{w}
\newcommand\FiltPt{z}

\begin{document}

\title{Knot Floer homology and the four-ball genus}

\author{Peter Ozsv\'ath\\Zolt{\'a}n Szab{\'o}}
\asciiauthors{Peter Ozsvath and Zoltan Szabo}
\coverauthors{Peter Ozsv\noexpand\'ath\\Zolt\noexpand\'an Szab\noexpand\'o}
\address{Department of
Mathematics, Columbia University\\New York 10025, USA}
\email{petero@math.columbia.edu}
\secondaddress{Department of
Mathematics, Princeton University\\New Jersey 08540, USA}
\secondemail{szabo@math.princeton.edu}

\asciiaddress{Department of
Mathematics, Columbia University\\New York 10025, USA\\and\\Department of
Mathematics, Princeton University\\New Jersey 08540, USA}
\asciiemail{petero@math.columbia.edu, szabo@math.princeton.edu}

\begin{abstract}  
We use the knot filtration on the Heegaard
Floer complex $\CFa$ to define an integer invariant $\tau(K)$ for
knots.  Like the classical signature, this invariant gives a
homomorphism from the knot concordance group to $\Z$. As such, it
gives lower bounds for the slice genus (and hence also the unknotting
number) of a knot; but unlike the signature, $\tau$ gives sharp bounds
on the four-ball genera of torus knots.  
As another illustration, we calculate the invariant for several
ten-crossing knots.
\end{abstract}
\asciiabstract{We use the knot filtration on the Heegaard Floer
complex to define an integer invariant tau(K) for knots.  Like the
classical signature, this invariant gives a homomorphism from the knot
concordance group to Z. As such, it gives lower bounds for the slice
genus (and hence also the unknotting number) of a knot; but unlike the
signature, tau gives sharp bounds on the four-ball genera of torus
knots.  As another illustration, we calculate the invariant for
several ten-crossing knots.}




\maketitlepage

\section{Introduction}

Let $K$ be a knot in the three-sphere. A slice surface for $K$ is a
smooth submanifold-with-boundary in $B^4$ whose boundary is the knot
$K$.  The four-ball genus $g^*(K)$ is the minimal genus of any slice
surface for $K$. The four-ball genus gives a lower bound for the
``unknotting number'' $u(K)$ (the minimal number of crossing-changes
required to unknot $K$).

Our aim here is to define an integer invariant of a knot $K$ which gives
a lower bound on $g^*(K)$, using the knot filtration
from~\cite{Knots}, see also~\cite{RasmussenThesis} and
\cite{Rasmussen}.  Specifically, letting $\CFa(S^3)$ denote the chain
complex whose homology calculates the Heegaard Floer homology
$\HFa(S^3)$ (cf.~\cite{HolDisk}, see also Section~\ref{sec:Review} for a
quick explanation), recall that a knot $K$ in $S^3$
induces a filtration on the chain complex $\CFa(S^3)$.  Letting
$\Filt(K,m)\subset \CFa(S^3)$ be the subcomplex generated by
intersection points whose filtration level is less than or equal to
$m$, we obtain an induced sequence of maps $$
\iota^m_K\colon H_*(\Filt(K,m))\longrightarrow
H_*(\CFa(S^3))=\HFa(S^3)\cong \Z,$$ which are isomorphisms for all
sufficiently large integers $m$. Let $$\tau(K)=\min\{m\in \Z
\big|\iota^m_K~{\text{~is non-trivial}}\}.$$

As we shall see, the invariant $\tau(K)$ in fact gives a lower bound
on the following generalization of the four-ball genus.  Let $W$ be a
smooth, oriented four-manifold with $\partial W = S^3$ and with
$b_2^+(W)=b_1(W)=0$ (as usual $b_1$ denotes the first Betti number, and
$b_2^+$ denotes the maximum dimension of any
vector subspace $V$ of $H^2(W)$ on which the cup-product form is
positive-definite). According to Donaldson's celebrated
theorem~\cite{DonKron}, the intersection form of $W$ is diagonalizable
(though in the applications it is interesting to consider the special
case where $W=B^4\#^b\mCP$, so the intersection form is already
diagonalized). Writing a homology class $[\Sigma]\in H_2(W)$ as
$$[\Sigma]=s_1\cm e_1+...+s_b \cm e_b,$$ where $e_i$ are an
ortho-normal basis for $H^2(W;\Z)$, and $s_i\in\Z$, we can define the
$L^1$ norm of $[\Sigma]$ by $$\Big|[\Sigma]\Big|=|s_1|+...+|s_b|.$$
Note that this is independent of the diagonalization (since the basis
$\{e_i\}$ is uniquely characterized, up to permutations and
multiplications by $\pm 1$, by the ortho-normality condition).  We
then have the following bounds on the genus of $\Sigma$, which are
proved in Section~\ref{sec:GenusBounds}:

\begin{theorem}
\label{thm:GenusBounds}
Let $W$ be a smooth, oriented four-manifold with $b_2^+(W)=0=b_1(W)$, 
and $\partial W = S^3$. If $\Sigma$ is any smoothly embedded
surface-with-boundary in $W$
whose boundary lies on $S^3$, where it is embedded as the knot $K$,
then we have the following inequality:
$$2\tau(K) +\Big|[\Sigma]\Big|+ [\Sigma]\cm [\Sigma]
\leq 2g(\Sigma).$$
\end{theorem}

The quantity $\tau(K)$ is additive under connected sums. This
additivity (together with the above theorem), can be reformulated in
the following terms.

Recall that two knots $K_1$ and $K_2$ are said to be concordant if
there is a smoothly embedded cylinder in $[1,2]\times S^3$ which meets
$\{i\}\times S^3$ in $K_i$ (for $i=1,2$). The connected sum of knots
descends to the set of concordance classes of knots $\Conc(S^3)$ to
endow the latter object with the structure of an Abelian group.  The
following result is also established in Section~\ref{sec:GenusBounds}:

\begin{theorem}
\label{thm:ConcordanceInvariance}
The map $\tau$ induces a group homomorphism from $\Conc(S^3)$ to $\Z$.
\end{theorem}

Indeed,  by letting $W$ be the four-ball in Theorem~\ref{thm:GenusBounds},
and reflecting $K$ if necessary, we obtain the following:

\begin{cor}
\label{cor:GenusBounds}
Let $K\subset S^3$ be a knot. Then,
$$|\tau(K)|\leq g^*(K).$$
\end{cor}

There is a classical knot invariant which has many of the same properties
which $\tau$ has --  the signature of $K$, $\sigma(K)$. In fact, $\tau$
and $-\sigma/2$ agree for a very wide class of knots. For instance, 
results from~\cite{AltKnots} (see also~\cite{Rasmussen}) give the following:

\begin{theorem}
\label{thm:SignatureTau}
Let $K$ be an alternating knot in $S^3$. Then, $\tau(K)=-\sigma(K)/2$.
\end{theorem}

\begin{proof}
From Theorem~\ref{AltKnots:thm:KnotHomology} of~\cite{AltKnots},
for an alternating knot, 
$H_*\left(\frac{\Filt(K,m))}{\Filt(K,m-1)}\right)$ is supported
in dimension $m+\frac{\sigma}{2}$. The result follows immediately.
\end{proof}

In fact, $\tau(K)=-\sigma(K)/2$ holds for some non-alternating knots as
well (compare Section~\ref{KT:sec:SmallKnots} of~\cite{calcKT}).

The following consequence of Theorems~\ref{thm:GenusBounds} and
\ref{thm:ConcordanceInvariance} proved in Section~\ref{sec:GenusBounds}
underscores the similarity between $-\sigma/2$ and $\tau$:

\begin{cor}
\label{cor:SkeinInequality}
Let $K_+$ be a knot in $S^3$, and $K_-$ be the new knot by changing
one positive crossing in $K_+$ to a negative crossing. Then,
$$\tau(K_+)-1\leq \tau(K_-)\leq \tau(K_+).$$
\end{cor}

Observe that $-\sigma(K)/2$ also satisfies the above
inequality. Indeed according to~\cite{Giller}, $-\frac{\sigma}{2}$ is
characterized among integer-valued knot invariants by the following
three properties:
\begin{itemize}
\item it vanishes for the unknot
\item it satisfies inequalities corresponding to the ones for $\tau$ stated
in the
above corollary
\item its parity is determined by the sign of
$\Delta_K(-1)$. \end{itemize}

However in general, $-\sigma/2\neq \tau$. As an example,
the knot $9_{42}$ has $\tau=0$ and $\sigma=2$.  (The
calculation of $\tau$ is a straightforward consequence of the
calculations from Proposition~\ref{Knots:prop:NineFortyTwo}
of~\cite{Knots}.) Other small examples are given in Section~\ref{sec:Examples}.
An infinite family of examples is provided 
by the following result proved in~\cite{NoteLens}:

\begin{theorem}
\label{thm:CalcTau}
[Corollary~\ref{NoteLens:cor:CalcTau} of~\cite{NoteLens}]
Let $K$ be a knot in $S^3$ and suppose that there is some integer $p\geq 0$ 
with the property that $S^3_p(K)=L(p,q)$ for some $q$. Then,
$\tau(K)$ is the degree of the symmetrized Alexander polynomial of $K$.
\end{theorem}

\begin{cor}
\label{cor:CalcTauTorusKnots}
Let $p$ and $q$ a pair of positive, relatively prime integers, and let
$T_{p,q}$ denote the $(p,q)$ torus knot. Then, 
$$\Big|\tau(T_{p,q})\Big|=\frac{pq-p-q+1}{2}.$$
\end{cor}

\begin{proof}
For a suitable chirality on $T_{p,q}$, we have that
$pq\pm 1$-surgery 
on $T_{p,q}$ is a lens space. Now, apply 
Theorem~\ref{thm:CalcTau}.
\end{proof}

Of course, the above corollary gives an infinite set of knots for
which $\sigma(K)\neq 2\tau(K)$ (see~\cite{GLM} for the calculation
of the signature of torus knots). For instance, $\sigma(T_{5,4})/2=4$,
while $\tau(T_{5,4})=5$.

It follows quickly from Corollary~\ref{cor:CalcTauTorusKnots} and
Corollary~\ref{cor:GenusBounds} that the four-ball genus and 
indeed the unknotting number of
$T_{p,q}$ is $(pq-p-q+1)/2$, a result first proved by Kronheimer and
Mrowka~\cite{KMMilnor} using Donaldson's invariants~\cite{DonKron},
and conjectured by Milnor (cf.~\cite{Milnor}, see
also~\cite{BoileauWeber}, \cite{Rudolph}).  Indeed, constructions of
Berge~\cite{Berge} give other fibered knots for which
Theorem~\ref{thm:CalcTau} applies.  Thus, for the knots arising from
Berge's constructions, the four-ball genus agrees with the degree of
the Alexander polynomial (which, since those knots are fibered, agrees
with their Seifert genus).  For more on this, see~\cite{NoteLens}.

We give also calculations for some ten-crossing knots in
Section~\ref{sec:Examples}. The calculations rest on the combinatorial
techniques developed in~\cite{AltKnots} and extended in~\cite{calcKT}.
These calculations can be used to determine the four-ball genera of
some $10$-crossing knots, whose four-ball genera were previously
calculated using techniques from gauge theory (cf.~\cite{KMMilnor},
\cite{Rudolph}, \cite{Tanaka},
\cite{Kawamura}).

In closing, we observe that many of the constructions from this paper
can be generalized to the case of null-homologous, oriented links in a
compact, oriented three-manifold. 
We sketch some of these generalizations in
Section~\ref{sec:Generalizations}.

In the proof of Theorem~\ref{thm:GenusBounds}, we make use of the knot
filtration and its relationship with the Heegaard Floer homology of a
corresponding surgered three-manifold. This relationship is spelled
out in detail in Section~\ref{Knots:sec:Relationship}
of~\cite{Knots}. For the reader's convenience, we review some aspects
of this in Section~\ref{sec:Review} before turning to the proofs of
the main result in Section~\ref{sec:GenusBounds}. The ten-crossing 
calculations are described in Section~\ref{sec:Examples}.

\vskip.2cm
\noindent{\bf{Remarks}}\qua 
The invariant $\tau(K)$ and also some additional constructions have
been independently discovered by Rasmussen in his
thesis~\cite{RasmussenThesis}.  A construction for obtaining
information on the slice genus, using Heegaard Floer homology in a
different way, has been developed by Strle and Owens,
cf.~\cite{OwensStrle}. We are also indebted to Jacob Rasmussen and
Chuck Livingston for many valuable comments on an early version of
this preprint.

\section{Properties of the knot filtration}
\label{sec:Review}

We recall here the knot filtration from~\cite{Knots}, focusing on the
case of knots in $S^3$. In Subsection~\ref{subsec:DefFilt}, we recall
the definition, and in Subsection~\ref{subsec:Surgeries} we recall its
relationship with Heegaard Floer homology of the surgered
three-manifold.

\subsection{Definition of the knot filtrations}
\label{subsec:DefFilt}

We briefly recall here the construction of the knot filtration.
For simplicity, we restrict to knots $K$ in the three-sphere.

Fix a doubly-pointed Heegaard diagram $(S,\alphas,\betas,w,z)$
for the knot $K\subset S^3$, in the following sense. Here, $S$ is
an oriented surface of genus $g$, $\alphas=\{\alpha_1,...,\alpha_g\}$
is a $g$ tuple of homologically linearly independent, pairwise
disjoint, simple closed curves in $S$, as is
$\betas=\{\beta_1,...\beta_g\}$. Of course, $\alphas$ and $\betas$
specify a pair of handlebodies $U_{\alpha}$ and $U_{\beta}$ which
bound $S$.  We require that $(S,\alphas,\betas)$ is a
Heegaard diagram for $S^3$, and also that the knot $K$ is supported
inside $U_\beta$ as an unknotted circle which meets the disk
attached to $\beta_1$ transversally in one point, and none of the other attaching
attaching disks. In particular, $\beta_1$ represents a meridian for $K$.  The
two points $w$ and $z$ lie on $S-\alpha_1-...-\alpha_g-\beta_1-...-\beta_g$, 
and can be connected by an arc
which meets $\beta_1$ exactly once, and none of the other attaching
circles. 

We consider the $g$-fold symmetric product $\Sym^g(S)$, with two
distinguished tori
\begin{eqnarray*}
\Ta=\alpha_1\times\cdots\times\alpha_g&{\text{and}}&
\Tb=\beta_1\times\cdots\times\beta_g.
\end{eqnarray*} The generators $X$ for the chain
complex $\CFa(S^3)$ are intersection points between $\Ta$ and $\Tb$
in $\Sym^g(S)$. 

Fix intersection points $\x, \y\in X$. A {\em Whitney disk} $u$
connecting $\x$ to $\y$
is a map 
$$u\colon \{\z\in\C\big| |z|\leq 1\} \longrightarrow \Sym^g(\Sigma)$$
satisfying the boundary conditions
\begin{eqnarray*}
u\{\zeta \big| \Real(\zeta)\geq 0~{\text{and}}~|\zeta|=1\}\subset \Ta,
&&
u\{\zeta\big| \Real(\zeta)\leq 0~{\text{and}}~|\zeta|=1\}\subset \Tb,\\
u(-\sqrt{-1})=\x, &&
u(\sqrt{-1})=\y.
\end{eqnarray*}
For a fixed point
$p\in S-\alpha_1-\cdots-\alpha_g-\beta_1-\cdots-\beta_g$, let $n_p(u)$
denote the algebraic intersection number of $u$ with the
submanifold $\{p\}\times \Sym^{g-1}(S)$. Note that $n_p(u)$ depends
only on the homotopy class $\phi$ of $u$. (In this context, homotopies
are to be understood as homotopies of maps all of which are
Whitney disks.)

There is a function $$\Filt \colon X \longrightarrow \Z$$ uniquely
characterized by the following two properties. For any $\x,\y\in
X$, we have that $$\Filt(\x)-\Filt(\y) =
n_\FiltPt(\phi)-n_\BasePt(\phi),$$ 
and also $$P(T)=\sum_{\x\in X}
\epsilon(\x)\cm T^{\Filt(\x)}$$ is a symmetric Laurent polynomial in the
formal variable $T$, where here, $\epsilon(\x)\in
\{\pm 1\}$ denotes the
local intersection number of $\Ta$ and $\Tb$ at $\x$, with respect to
fixed orientations on the tori and $\Sym^g(S)$.
(Indeed, for the appropriate choice of overall sign, $P(T)$ coincides
with the symmetrized Alexander polynomial of $K$, cf.\
Equation~\eqref{Knots:eq:EulerChar} of~\cite{Knots}.)

Recall~\cite{HolDisk} that there is a homology theory for (closed,
oriented) three-manifolds $\HFa(Y)$, whose generators are intersection
points $X$, and whose differentials count pseudo-holomorphic Whitney
disks in $\Sym^g(S)$, in the homotopy class with
$n_\BasePt(\phi)=0$. In the case where $Y\cong S^3$, $\HFa(Y)\cong
\Z$.  Starting with a doubly-pointed 
Heegaard diagram for $S^3$ compatible with a
knot $K\subset S^3$, it is not  difficult to see
that if we let $\Filt(K,m)\subset \CFa(S^3)$ denote the subset
generated by those $\x\in X$ with $\Filt(\x)\leq m$, then that subset
indeed is preserved by the differential, i.e.\ $\Filt(K,m)$ defines a
$\Z$-filtration of $\CFa(S^3)$, indexed by integers $m\in\Z$.
We let $\CFKa(S^3,K)$ denote the chain complex of $\CFa(S^3)$,
together with this $\Z$-filtration induced from $K$.

Since there are only finitely many generators in $X$, the filtration
we have defined has finite support; i.e.\ for all sufficiently small
$m\in \Z$, $\Filt(K,m)=0$, and for all sufficiently large $m\in \Z$,
$\Filt(K,m)=\CFa(S^3)$, and in particular, if $\iota_K^m$ denotes the
map on homology induced by the inclusion $\Filt(K,m)$ in $\CFa(S^3)$,
then $\iota_K^m$ is trivial for all sufficiently small $m$, and an
isomorphism for all sufficiently large $m$, and hence the quantity
$\tau(K)$ defined in the introduction is a finite integer.  Although
$\tau(K)$ as defined might appear to depend on a choice of particular
Heegaard diagram used to define the knot filtration, it is shown in
Theorem~\ref{Knots:thm:KnotInvariant} of~\cite{Knots} that in fact the
filtered chain homotopy type of the filtered complex $\CFKa(S^3,K)$ is a
knot invariant, and and hence so is the integer $\tau(K)$. (Actually,
Theorem~\ref{Knots:thm:KnotInvariant} of~\cite{Knots} is explicitly
stated for the induced filtration of $\CFinf$, which generalizes the
filtration of $\CFa$ we just described, see also the discussion
below.)

Recall that $\HFa(Y)$ is the simplest of the Heegaard Floer homologies
associated to three-manifolds in~\cite{HolDisk}. There are also
induced knot filtrations on the chain complexes associated to the
other variants of Heegaard Floer homology.  Although these
filtrations are not used in the definition of $\tau$, they are used in
the verification of its four-dimensional properties.

To this end, recall that in~\cite{HolDisk}, there is another invariant
for integer homology three-spheres $Y$, 
$\CFinf(Y)$, whose generators are pairs
$[\x,i]\in\left(\Ta\cap\Tb\right)\times \Z$, endowed with a differential
given by the formula
$$\partial [\x,i] = \sum_{\y\in\Ta\cap\Tb} \sum_{\{\phi\in\pi_2(\x,\y)\}}
\#\left(\frac{\ModFlow(\phi)}{\R}\right) [\y,i-n_\BasePt(\phi)],$$
where here $\pi_2(\x,\y)$ denotes the space of homotopy classes of
Whitney disks in $\Sym^g(S)$ connecting $\x$ to $\y$, and
$\ModFlow(\phi)$ denotes the moduli space of pseudo-holomorphic
representatives for $\phi$, and in fact, the coefficient of
$[\y,i-n_\BasePt(\phi)]$ is a suitable signed count of points in the
moduli space, after we divide out by the one-dimensional 
automorphism
group of the unit disk in the complex plane, fixing
$\sqrt{-1}$ and $-\sqrt{-1}$. This chain complex admits a subcomplex $\CFm(Y)$, generated
by those pairs $[\x,i]$ with $i<0$, and a quotient complex $\CFp(Y)$,
which can be thought of as generated by pairs $[\x,i]$ with $i\geq 0$.
There is an endomorphism $U$ of the chain complex $\CFinf(Y)$, which
respects the subcomplex $\CFm(Y)$, defined by $U[\x,i]=[\x,i-1]$.

Indeed, the same definition can be made for rational homology spheres.
In this case, there is an identification $\spinc_\BasePt\colon
X\longrightarrow \SpinC(Y)$, and a corresponding splitting of
complexes $$\CFinf(Y)\cong
\bigoplus_{\spinc\in\SpinC(Y)}\CFinf(Y,\spinc).$$ Extra care is to be
taken in the case where $b_1(Y)>0$, but we do not describe this here,
as such three-manifolds do not play a role in the present article.
The main result from~\cite{HolDisk} states that the homology groups
$\HFm(Y,\spinc)$, $\HFinf(Y,\spinc)$, and $\HFp(Y,\spinc)$, thought of
as modules over $\Z[U]$, are topological invariants of the three-manifold
$Y$ and its $\SpinC$ structure $\spinc$.
In fact,  the long exact sequences associated to the two
short exact sequences
\begin{equation}
\begin{CD}
\label{eq:CanonicalSequences}
0@>>>\CFm(Y,\spinc) @>>>\CFinf(Y,\spinc) @>>>\CFp(Y,\spinc) @>>>0  \\
0@>>>\CFa(Y,\spinc) @>>>\CFp(Y,\spinc) @>{U}>>\CFp(Y,\spinc) @>>>0  \\
\end{CD}
\end{equation}
are also topological invariants.

A knot $K\subset S^3$ induces a $\Z\oplus\Z$ filtration of $\CFinf(S^3)$
(and also $\CFm(S^3)$ and $\CFp(S^3)$), as follows.
Let $\CFKinf(S^3,K)$ denote the complex generated by
$[\x,i,j]\in\left(\Ta\cap\Tb\right) \times \Z\times \Z$ with
the property that $\Filt(\x)+(i-j)=0$, and differential $$\partial
[\x,i,j] =
\sum_{\y\in\Ta\cap\Tb} \sum_{\{\phi\in\pi_2(\x,\y)\}}
\#\left(\frac{\ModFlow(\phi)}{\R}\right) [\y,i-n_\BasePt(\phi),j-n_\FiltPt(\phi)]$$
(where, once again, we use some doubly-pointed Heegaard
diagram compatible with the knot).  The map which which associates to
$[\x,i,j]$ the pair $(i,j)\in\Z\oplus \Z$ induces a
$\Z\oplus\Z$-filtration on $\CFKinf(S^3,K)$, meaning that if
$[\y,k,\ell]$ appears in $\partial [\x,i,j]$ with a non-zero
coefficient, then  $k\leq i$ and $\ell\leq j$.  There is an
endomorphism of $\CFKinf(S^3,K)$ defined by $U\cm [\x,i,j]=[\x,i-1,j-1]$.
The forgetful map sending $[\x,i,j]$ to $[\x,i]$ induces an
isomorphism of chain complexes $\CFKinf(S^3,K)$ and $\CFinf(S^3)$.
Moreover, this map is $\Z[U]$-equivariant.

\subsection{Knot filtrations and surgeries}
\label{subsec:Surgeries}

Let $K\subset Y$ be a knot with framing $\lambda$ in a three-manifold,
and let $X_\lambda(K)$ denote the four-manifold obtained by attaching a
two-handle to $[0,1]\times Y$ to $\{1\}\times Y$ along $K$ with
framing $\lambda$. This can be thought of as a cobordism from $Y$ to
the three-manifold $Y_\lambda(K)$ obtained by performing
$\lambda$-framed surgery on $Y$ along $K$. Given an $\SpinC$ structure
$\spinc$ on $X_\lambda(K)$, there is an induced map $${\widehat
F}_{W,\spinc}\colon \HFa(Y,\spinc|_Y) \longrightarrow
\HFa(Y_{\lambda}(K),\spinc|_{Y_\lambda(K)}),$$
(and also maps for the other versions of Floer homology). This map is
induced from a corresponding chain map, gotten by counting
pseudo-holomorphic triangles in $\Sym^g(\Sigma)$, as explained in
Section~\ref{HolDiskTwo:sec:Surgeries} of~\cite{HolDiskTwo}. Further
invariance properties of these maps, and a generalization to other
cobordisms, are established in~\cite{HolDiskFour}.

In Section~\ref{Knots:sec:Relationship} of~\cite{Knots}, we described
the relationship between this knot filtration and the Heegaard Floer
homologies of three-manifolds obtained by performing ``sufficiently
large'' integral surgeries on $S^3$ along $K$. We sketch the results
here, and refer the reader to~\cite{Knots} for a more thorough
treatment. Moreover, this relationship gives an interpretation of the
some of maps induced by cobordisms in terms of the knot filtration.

Let $K\subset S^3$ be a knot. We write simply $C$ for the $\Z\oplus\Z$
filtered complex $\CFKinf(S^3,K)$, suppressing the knot from the
notation.  The subcomplex $\CFm(S^3)$ is represented by the complex
$C\{i< 0\}$. (This notation means that the complex in question
is the subcomplex of $C$ generated by those homogeneous elements whose
$\Z\oplus\Z$ filtration level $(i,j)$ satisfies the stated
constraint). Its quotient $\CFp(S^3)$ is represented by $C\{i\geq
0\}$. The subcomplex $\CFa(S^3)\subset\CFp(S^3)$ is represented by
$C\{i=0\}$, and the various filtration levels $\Filt(K,m)$ correspond
to the subcomplexes $C\{i=0,j\leq m\}$ of $C\{i=0\}$ (i.e.\ the
subcomplex generated $[\x,i,j]$ where with $i=0$ and $j\leq m$).

Framings on knots in $S^3$ are canonically identified with the
integers. 
For given $n\in\Z$, on the three-manifold $S^3_{-n}(K)$
obtained by $-n$-framed surgery on $S^3$ along $K$, there is a natural
affine identification $\SpinC(S^3_{-n}(K))\cong\Zmod{n}$, specified
by an orientation for the knot $K$. More
precisely, an orientation of $K$ induces an orientation on its Seifert
surface $\Sigma$. The oriented Seifert surface $\Sigma$ can be capped
off inside the two-handle to obtain a closed surface ${\widehat
\Sigma}$ inside the cobordism $X_{-n}(K)$ from $S^3$ to $S^3_{-n}(K)$
obtained by the two-handle addition. Now, for fixed $m\in\Z$, consider
the $\SpinC$ structure $\spinc_m\in\SpinC(X_{-n}(K))$ with the property
that $$\langle c_1(\spinc_m),[{\widehat \Sigma}]\rangle - n = 2m.$$
Restricting to $S^3_{-n}(K)\subset \partial X_{-n}(K)$, we get the induced identification
$$\SpinC(S^3_{-n}(K))\cong
\Zmod{n}.$$
For $[m]\in\Zmod{n}$, we let $\HFa(S^3_{-n}(K),[m])$ denote the summand
of the Floer homology corresponding to the $\SpinC$ structure
corresponding to $[m]$.  (Although the orientation on $K$ might appear
important here, the other choice of orientation induces a conjugate
identification of $\SpinC$ structures, which gives an isomorphic theory, cf.\
Proposition~\ref{Knots:prop:JInvarianceGen} of~\cite{Knots}.)

Theorem~\ref{Knots:thm:LargeNegSurgeries} of~\cite{Knots}
states that for each integer
$m\in\Z$, there is an integer $N$ so that for all $n\geq N$,
the short exact sequence of complexes
\begin{equation}
\label{eq:StandardSequence}
0\to \CFm(S^3_{-n}(K),[m]) \longrightarrow
\CFinf(S^3_{-n}(K),[m]) \longrightarrow 
\CFp(S^3_{-n}(K),[m]) \to 0
\end{equation}
is identified with the natural short exact sequence
$$
0\longrightarrow C\{\min(i,j-m)<0\} \longrightarrow C \longrightarrow C\{\min(i,j-m)\geq 0 \} \longrightarrow0
.$$
Similarly, the inclusion
\begin{equation}
\label{eq:StandardInclusion}
\begin{CD}
0@>>>\CFa(S^3_{-n}(K),[m]) @>>> \CFp(S^3_{-n}(K),[m])
\end{CD}
\end{equation}
is identified with the inclusion
$$
\begin{CD}
0@>>>C\{\min(i,j-m)=0\} @>>>C\{\min(i,j-m)\geq 0\}.
\end{CD}
$$
(In the interest of simplicity, we
will remain in the context of $\HFa$ as much as possible, and in fact
will not to use Theorem~\ref{Knots:thm:LargeNegSurgeries} for
the theories $\HFm$, $\HFinf$, and $\HFp$.)

There is a natural chain map $C\{i=0\}\longrightarrow
C\{\min(i,j-m)=0\}$, which is defined to vanish in the subcomplex of
$C\{i=0,j<m\}$. In fact, the
proof of Theorem~\ref{Knots:thm:LargeNegSurgeries} of~\cite{Knots}
shows that this chain map induces the map from $\HFa(S^3)$
to $\HFa(S^3_{-n}(K),[m])$ given by the two-handle addition, endowed with the
$\SpinC$ structure for which $$\langle c_1(\spinc),[{\widehat
\Sigma}]\rangle-n=2m$$ 
(again, provided that $n$ is sufficiently large compared to $m$ and the
genus of the knot).

Note that a similar picture holds for positive integral surgeries
(cf.\ Theorem~\ref{Knots:thm:LargeNegSurgeries}), but we do not need
that statement here.

\section{Proof of Theorem~\ref{thm:GenusBounds}}
\label{sec:GenusBounds}

Our aim here is to prove Theorem~\ref{thm:GenusBounds}. As a
preliminary step, we establish a four-dimensional interpretation of
$\tau$ which will be useful later.  Then, we establish some {\em a
priori} properties of $\tau$: additivity under connected sums, and its
behavior under reflection. As a third step, we describe some
properties of the maps induced on $\HFa$ by cobordisms (analogous to
the usual adjunction inequalities of four-manifold invariants).  After
establishing these preliminaries, we prove
Theorem~\ref{thm:GenusBounds}. Proofs of its corollaries are given in
the end of the section.

\subsection{A four-dimensional interpretation of $\tau$}
The following interpretation of $\tau$ 
will be useful to us.

But first, we set up notation.
For $m\in \Z$, we have a short exact sequence
$$\begin{CD}
0@>>>\Filt(K,m)@>{I^m_K}>>\CFa(S^3) @>{P^m_K}>>Q(K,m)@>>>0,
\end{CD}$$
where $I^m_K$ is the natural inclusion, and $P^m_K$ and $Q(K,m)$ are 
are defined to make the sequence exact. Let $\iota^m_K$ and $p^m_K$
denote the maps induced by $I^m_K$ and $P^m_K$ on homology. Of course,
$I^m_K$ is represented by the inclusion of $C\{i=0,j\leq m\}$ inside
$C\{i=0\}$.

Fix a knot $K\subset S^3$, and an integer $n$.
Let
$${\widehat F}_{n,m}\colon \HFa(S^3)\longrightarrow \HFa(S^3_{-n}(K),[m])$$
denote the map associated to the two-handle addition,
endowed with the $\SpinC$ structure $\spinc_m$ characterized by
$$\langle c_1(\spinc_m),[{\widehat \Sigma}]\rangle-n=2m.$$

\begin{prop}
\label{prop:FourDInterp}
If $m<\tau(K)$, then ${\widehat F}_{n,m}$ is non-trivial for
all sufficiently large $n$. Also, if $m>\tau(K)$, then ${\widehat F}_{n,m}$
is trivial for all sufficiently large $n$.
\end{prop}

\begin{proof}
Consider the diagram
\begin{equation}
\label{eq:Diag1}
\begin{CD}
0\to \Filt(K,m) @>{I^m_K}>> C\{i=0\}\simeq\CFa(S^3) @>{P^m_K}>> Q(K,m)\to0\\
 @V{\Pi}VV @V{f}VV @V{\cong}VV \\ 
0\to C\{i\geq 0,j=m\} @>>> C\{\min(i,j-m)=0\} @>>>Q(K,m) \to0. 
\end{CD}
\end{equation}
According to Theorem~\ref{Knots:thm:LargeNegSurgeries} of~\cite{Knots}
for $n$ sufficiently large, we have an identification
$$C\{\min(i,j-m)=0\}\simeq\CFa(S^3_{-n}(K),[m])$$ under which the map
$f$ represents the chain map ${\widehat F}_{n,m}$ above.  Now, if
$m<\tau(K)$, the induced map on homology map $p^m_K$ induces an
injection in homology, and hence ${\widehat F}_{m,n}$ is non-trivial.
Moreover, since the projection $f$ is trivial on $$C\{i=0,j\leq
m-1\}=\Filt(K,m-1),$$ ${\widehat F}_{n,m}$ factors through the map $P^{m-1}_K$.
If $m>\tau(K)$, the induced map $p^{m-1}_K$ on homology is trivial,
and hence so is ${\widehat F}_{n,m}$.
\end{proof}

\subsection{Additivity of $\tau$}
The additivity of $\tau$ under connected sums follows readily from the
K\"unneth principle for the knot filtration,
Theorem~\ref{Knots:thm:ConnectedSumsOfKnots} of~\cite{Knots}.

The maps $I^m_K$ induces a filtration of $\CFa(S^3)\otimes_\Z
\CFa(S^3)$ as the image of $$
\sum_{m_1+m_2=m}I^{m_1}_{K_1}\otimes I^{m_2}_{K_2}\colon \!\!\!
\bigoplus_{m_1+m_2=m} \!\!\!\Filt(K_1,m_1)\otimes_\Z
\Filt(K_2,m_2) \longrightarrow \CFa(S^3)\otimes_\Z\CFa(S^3).$$ 

According to Theorem~\ref{Knots:thm:ConnectedSumsOfKnots}
of~\cite{Knots}, under a homotopy equivalence
$\CFa(S^3)\otimes_\Z\CFa(S^3)\simeq\CFa(S^3)$, the above filtration is
identified with the filtration of $\CFa(S^3)$ induced by the connected
sum $K_1\# K_2$.

Indeed, for simplicity of exposition, we switch from base ring $\Z$ to
$\Q$, a change which we also suppress from the notation. In this case, then,
we can think of $\tau(K)$ as the minimum integer for which 
$$\iota^m_K\colon H_*(\Filt(K,m))\longrightarrow \HFa(S^3)\cong \Q$$
is surjective.

\begin{prop}
\label{prop:Additivity}
Let $K_1$ and $K_2$ be a pair of knots in $S^3$, and let $K_1\# K_2$
denote their connected sum. Then,
$$\tau(K_1\# K_2)=\tau(K_1)+\tau(K_2).$$
\end{prop}

\begin{proof}
According to the theorem quoted above,
$I^m_{K_1\#K_2}$ is surjective on homology if and only if there
is a decomposition $m=m_1+m_2$ with the property that
$$f=(I^{m_1}_{K_1}\otimes_\Q I^{m_2}_{K_2})_*
\colon H_*(\Filt(K_1,m_1)\otimes_\Q \Filt(K_2,m_2)) \longrightarrow
\HFa(S^3)\cong \Q$$
is surjective. Now, by the K\"unneth formula, we have
an identification
$$H_*(\Filt(K_1,m_1))\otimes_\Q H_*(\Filt(K_2,m_2))\cong
H_*(\Filt(K_1,m_1)\otimes_\Q \Filt(K_2,m_2)),$$
and hence the map $f$ is surjective if and only if
\begin{align*}
\iota^{m_1}_{K_1}\otimes_\Q \iota^{m_2}_{K_2}\colon H_*(\Filt(K_1,m_1))
\otimes_\Q
H_*(\Filt(K_2,m_2))
&\longrightarrow \HFa(S^3)\otimes_\Q\HFa(S^3)\\
&\hspace{2cm}\cong \HFa(S^3)\cong \Q
\end{align*}
is. This in turn is easily seen to be surjective if and only if both
$\iota^{m_1}_{K_1}$ and $\iota^{m_2}_{K_2}$ are surjective.
This shows that $\tau(K_1\#
K_2)=\tau(K_1)+\tau(K_2)$. 
\end{proof}

It is worth pointing out that, strictly speaking, the knot filtration
depends on an oriented knot. However, conjugation invariance
of the knot filtration shows that $\tau(K)$ is independent of this 
additional choice (cf.\ Proposition~\ref{Knots:prop:OrientKnot} of~\cite{Knots}).

We have also the following result:

\begin{lemma}
\label{lemma:Reflection}
Let $K$ be a knot, and let $-K$ denote its reflection. Then, 
$$\tau(-K)=-\tau(K).$$
\end{lemma}

\begin{proof}
Let $I^m_K$ be as before, and indeed we have a short exact sequence
$$
\begin{CD}
0@>>> \Filt(K,m) @>{I^m_K}>> \CFa(S^3) @>{P^m_K}>> Q(K,m) @>>> 0.
\end{CD}
$$

According to~\cite{HolDiskTwo}, there is a duality map $$\Duality\colon
\HFa_*(S^3)\longrightarrow \HFa^*(-S^3),$$ induced by a map of chain
complexes, which we also denote by $$\Duality\colon
\CFa_*(S^3)\longrightarrow \CFa^*(-S^3).$$ It is easy to see that
under this map, if $K\subset S^3$ is a knot, then we have the commutative
diagram
$$
\begin{CD}
\Filt(K,m) @>{I^m_K}>> \CFa_*(S^3) \\
@V{\Duality}V{\cong}V  @V{\cong}V{\Duality}V \\
Q^*(-K,-m) @>{P^{-K}_{-m}}>> \CFa^*(S^3), \\
\end{CD}
$$
where
$P^{-m}_{-K}$ is the map which is dual to $P^{-m}_{-K}$. 
The induced map on cohomology $p^{-m}_{-K}$ is trivial if and only if
othe map $\iota^{-m}_{-K}$ is non-trivial. The lemma now follows.
\end{proof}

\subsection{Maps on $\HFa$}
We turn now to some more lemmas which will be used in the proof of
Theorem~\ref{thm:GenusBounds}. But first, we must set up some more notation.

Let $W$ be a four-manifold with $b_2^+(W)=0=b_1(W)$ and $\partial W=S^3$, and
let $K\subset S^3$ be a knot. Let $W_{-n}(K)$ denote the four-manifold
obtained by attaching a two-handle to $W$ along $K$, with framing
$-n$. As in the statement of Theorem~\ref{thm:GenusBounds}, we fix a 
surface $\Sigma$ whose boundary lies in $\partial W$, where it 
agrees with the knot $K$. This surface-with-boundary can be closed off
to obtain a smoothly embedded surface ${\widehat \Sigma}$ inside
$W_{-n}(K)$. 

Then, we can view $W-B^4$ as a cobordism from $S^3$ to
$S^3$. According to Donaldson's theorem, $W$ has diagonalizable
intersection form. Thus, if we let $b=b_2(W)$, there are are $2^b$
characteristic vectors $K$ for the intersection form (on 
$H^2(W;\Z)/\Tors$) with $K\cm K =
-b$. Note that if $[\Seif]\in H_2(W;\Z)$, then $$\Big| [\Seif]\Big| =
\max_{\{K\in\Char(W)\big| K\cm K = -b\}}
\langle K, [\Seif]\rangle,$$
where $\Char(W)\subset H^2(W;\Z)/\Tors$ denotes the set of characteristic 
vectors for the intersection form.

Note that in $W_{-n}(K)$,
$$[{\widehat \Sigma}]\cm[{\widehat \Sigma}] = 
[{\Sigma}]\cm [{\Sigma}]
-n.$$

\begin{lemma}
\label{lemma:NegDef}
Let $\spinc$ be a $\SpinC$ structure over a four-manifold $W$ 
with $\partial W = S^3$, $b_2^+(W)=0$, and $b_2(W)=b$.
Then its first Chern class satisfies
\begin{equation}
\label{eq:MinimalClass}
c_1(\spinc)\cm c_1(\spinc)=-b
\end{equation}
if and only if the induced map
$${\widehat F}_{W-B,\spinc}\colon \HFa(S^3)\longrightarrow \HFa(S^3)$$
is non-trivial (in which case it is an isomorphism).
\end{lemma}

\begin{proof}
In~\cite{AbsGraded} (see especially the proof of
Theorem~\ref{AbsGraded:thm:IntFormQSphere} of~\cite{AbsGraded}), it is
shown that the map induced by $W-B^4$ on $\HFinf$ is an isomorphism,
and its shift in degree is given by $$(c_1(\spinc)^2+b)/4.$$ Thus, the
map on $\HFp$ induced by a $\SpinC$ structure $\spinc$ is an
isomorphism if and only if Equation~\eqref{eq:MinimalClass} holds. The
lemma now follows readily from the long exact sequence relating $\HFa$
and $\HFp$ (cf.\ Equation~\eqref{eq:CanonicalSequences}), and its
functoriality under the maps induced by cobordisms (cf.~\cite{HolDiskFour}).
\end{proof}

\begin{lemma}
\label{lemma:TubularNbd}
Let $N$ be the total space of circle bundle with Euler number $-n<0$
over an oriented
two-manifold $\Sigma$ of genus $g>0$.
The map
$${\widehat F}_{N-B,\spinc}\colon \HFa(S^3)\longrightarrow \HFa(\partial N)$$
is trivial whenever
$$\langle c_1(\spinc), [\Sigma]\rangle+[\Sigma]\cm[\Sigma] > 2g(\Sigma)-2.$$
\end{lemma}

\begin{proof}
We argue as in~\cite{HolDiskSymp} and~\cite{SympThom}, making use
of the absolute grading on Floer homology. Specifically,
according to Theorem~\ref{HolDiskFour:thm:AbsGrade}
of~\cite{HolDiskFour}, there is an absolute $\Q$-grading 
on the Floer homology 
$\HFa$  of any three-manifold equipped with a torsion $\SpinC$
structure, which is uniquely characterized by the 
following two properties: $\HFa(S^3)$ is supported in dimension zero,
and if $W$ is a cobordism 
from $Y_1$ to $Y_2$, which is given a $\SpinC$ structure $\spinc$ whose restrictions to $Y_1$ and $Y_2$ have torsion first Chern class, then 
the induced map
$${\widehat F}_{W,\spinc}\colon \HFa(Y_1,\spinc|_{Y_1})\longrightarrow
\HFa(Y_2,\spinc|_{Y_2})$$
shifts degree by
\begin{equation}
\label{eq:DegreeShift}
\frac{c_1(\spinc)^2-2\chi(W)-3\sigma(W)}{4}.
\end{equation} With
respect to this absolute grading, the rank of
$\HFa_{i}(\#^{2g}(S^2\times S^1))$ is zero if $|i|> g$. 
The Heegaard Floer homology of 
$\#^{2g}(S^2\times S^1)$ can be calculated directly from its Heegaard diagram,
cf.\
Subsection~\ref{HolDisk:subsec:STwoTimesSOne} of~\cite{HolDisk}. 
Indeed there
is a constant $c$ with the property that $\HFa_{i+c}(\#^{2g}(S^2\times
S^1))$ is a free module whose rank is given by the binomial
coefficient $\left(\begin{array}{c} 2g \\ i\end{array}\right)$.
(We suppress $\SpinC$ structures
from the notation for $\HFa$ of $\#^{2g}(S^2\times S^1)$,
since that is non-trivial only in the $\SpinC$
structure with trivial first Chern class;
the subscript $i$ here denotes the absolute $\Q$-degree). 
The
stated vanishing follows from the fact that $c=g$. To see this, note
that a direct inspection of the Heegaard diagrams
shows that  the non-zero elements in $\HFa(\#^{2g}(S^2\times
S^1))$ with lowest degree are in the image of the map on $\HFa$
induced by the
cobordism $S^3$ to $\#^{2g}(S^2\times S^1)$ obtained by attaching
${2g}$ two-handles. In turn this cobordism shifts degrees down by $g$,
according to Equation~\eqref{eq:DegreeShift}.

Now, by blowing up $N$ sufficiently many times (and using the
``blow-up formula'' of~\cite{HolDiskFour}), we reduce to the case
where $N$ is the total space of a circle bundle with $-n<-2g+1<0$. In
this case,
\begin{equation}
\label{eq:FloerHomologyCircleBundle}
\HFa(\partial N,\spinc|\partial N)\cong \HFa(\#^{2g}(S^2\times
S^1)),
\end{equation}
as relatively $\Z$-graded Abelian groups, as can be seen by appealing
to the long exact sequence for integral surgeries. In fact, letting
$\spinc_0$ denote the $\SpinC$ structure which minimizes
$c_1(\spinc_0)^2$ among all $\SpinC$ structures $\spinc'$ over $N$
which have $\spinc'|_{\partial N}=\spinc|_{\partial N}$, it is the map
induced by $\spinc_0$, which shifts absolute degree by
$\frac{c_1(\spinc_0)^2+1}{4}$, which induces the isomorphism of
Equation~\eqref{eq:FloerHomologyCircleBundle}:  
$${\widehat F}_{W,\spinc_0|W}\colon
\HFa(\#^{2g}(S^2\times S^1))\stackrel{\cong}{\longrightarrow}
\HFa(\partial N,\spinc|\partial N).$$
Here we have broken $N$ into one zero handle, $2g$ one-handles, and one
two-handle; it is the latter which specifies the cobordism
$W$ from $\#^{2g}(S^2\times S^1)$ to $\partial N$.
(Essentially this calculation, with more details,
can be found in Lemma \ref{AbsGraded:lemma:CorrTermCircleBundle}
of~\cite{AbsGraded}, only there we consider the case of $\HFp$, 
rather than $\HFa$.)

The stated hypothesis on $\spinc$, together with the fact that
$\HFa(\#^{2g}(S^2\times S^1))$ is supported in degrees $[-g,g]$ shows now
that the map induced by $\spinc$ is trivial.
\end{proof}

\subsection{Proof of Theorem~\ref{thm:GenusBounds}}

Let $W$ and $\Sigma$ be as in the statement of the theorem.  
We subdivide the proof into two cases: $g(\Sigma)=g>0$ and $g=0$.

\vskip.3cm
\noindent{\bf{Proof of Theorem~\ref{thm:GenusBounds} when $g>0$}}\qua
For an integer $n$ (which we will fix later), we let $W'$ denote the
four-manifold obtained by deleting an open four-ball from the interior
of four-manifold $W_{-n}(K)$ (which is disjoint from $\Sigma$).  This
four-manifold decomposes as $$W'\cong W_1\cup_{S^3} W_2,$$ where $W_1$
is obtained by deleting a small four-ball from $W$ (and hence it is
independent of $n$), while $W_2=X_{-n}(K)$ is the cobordism from $S^3$
to $S^3_{-n}(K)$ specified by the two-handle addition. We close off
$\Sigma$ inside the two-handle to obtain an surface ${\widehat
\Sigma}$ with $g({\widehat \Sigma})=g(\Sigma)$. Moreover, we
can split 
the homology class $[{\widehat \Sigma}]=[\Sigma_1]\oplus[\Sigma_2]$,
where $[\Sigma_i]\in H_2(W_i;\Z)$.
Note that 
\begin{eqnarray*}
\Big|[\Sigma_1]\Big|=\Big|[\Sigma]\Big| &{\text{and}}&
[\Sigma]\cm[\Sigma]=[\Sigma_1]\cm [\Sigma_1].
\end{eqnarray*}

Fix a $\SpinC$ structure $\spinc_1\in\SpinC(W_1)$ so that
\begin{eqnarray*}
c_1(\spinc_1)^2+b=0 &{\text{and}}&
\Big|[\Sigma_1]\Big|=\langle c_1(\spinc_1),[\Sigma_1]\rangle.
\end{eqnarray*}
For any integer $m<\tau(K)$, we choose $n$ 
large enough that Proposition~\ref{prop:FourDInterp}
holds, and fix a
$\SpinC$ structure $\spinc_2$ over $W_2$ so that $$\langle
c_1(\spinc_2),[\Sigma_2]\rangle-n = 2m<2\tau(K).$$ According to
Lemma~\ref{lemma:NegDef} and Proposition~\ref{prop:FourDInterp} 
respectively, the
maps ${\widehat F}_{W_1,\spinc_1}$ and ${\widehat F}_{W_2,\spinc_2}$
induce non-trivial maps on $\HFa$. Thus, by naturality of
the maps induced by cobordisms
(cf.~\cite{HolDisk}), if we let $\spinc$ be a $\SpinC$ structure
with $\spinc|W_i=\spinc_i$,  the map ${\widehat F}_{W',\spinc}$ 
induces a non-trivial map on $\HFa$. 

Note that  ${\widehat \Sigma}$ is represented by a closed,
embedded surface of
genus $g$, so we can split the cobordism $W'$ alternately as $W_1'\cup
W_2'$, where $W_1'$ is the tubular neighborhood of ${\widehat \Sigma}$
minus a four-ball
(thought of as a subset of $W'$, containing one of its boundary
components), and $W_2'$ is the remaining part of $W'$.

It now follows from Lemma~\ref{lemma:TubularNbd} that 
$$
\langle c_1(\spinc),[{\widehat \Sigma}]\rangle + [{\widehat \Sigma}]\cm
[{\widehat \Sigma}] \leq 2g-2.$$
On the other hand, it is easy to see that the left-hand-side is
$$\Big|[\Sigma_1]\Big| + [\Sigma_1]\cm[\Sigma_1] 
+ \langle c_1(\spinc),[\Sigma_2]\rangle + [\Sigma_2]\cm[\Sigma_2]
= \Big|[\Sigma_1]\Big| + [\Sigma_1]\cm[\Sigma_1] + 2m.$$

It follows at once that if $g>0$, then
$$ 2\tau(K)+\Big|[\Sigma]\Big|+[\Sigma]\cm[\Sigma]\leq g(\Sigma).\eqno{\qed} $$

\vskip.3cm
\noindent{\bf{Proof of Theorem~\ref{thm:GenusBounds} when $g=0$}}\qua
Fixing $W$, $K$, and $\Sigma$ as before, 
except now we assume that $\Sigma$ is a disk.
We can form a new four-manifold
$W\#_b W$ by boundary connected sum, which contains $K\# K$ on its boundary,
a knot which bounds the smoothly-embedded disk $\Sigma'=\Sigma\#_b\Sigma$.
Adding a trivial handle to $\Sigma'$, and using the previous case
of the theorem, we see that
$$2\tau(K\# K) + \Big|[\Sigma']\Big| + [\Sigma']\cm[\Sigma']\leq 2.$$
According to Proposition~\ref{prop:Additivity},
$\tau(K\# K)=2\tau(K)$, and it is also
easy to see that 
\begin{eqnarray*}
\Big|[\Sigma']\Big|=2\Big|[\Sigma]\Big| &{\text{and}}&
[\Sigma']\cm[\Sigma']=2[\Sigma]\cm[\Sigma].
\end{eqnarray*}
Thus, we see that
$$2\tau(K)+\Big|[\Sigma]\Big|+[\Sigma]\cm[\Sigma]\leq 1.$$
But the left-hand-side is easily seen to be an even integer,
so the stated inequality immediately follows.
\qed

\subsection{Corollaries.}

We now turn to some consequences of Theorem~\ref{thm:GenusBounds}.

\vskip.2cm
\noindent{\bf{Proof of Corollary~\ref{cor:GenusBounds}}}\qua
Apply Theorem~\ref{thm:GenusBounds} with $W=B^4$, so that $[\Sigma]=0$,
to see that $\tau(K)\leq g^*(K)$. Reflecting the knot $K$ and applying
Lemma~\ref{lemma:Reflection}, we get that $-\tau(K)\leq g^*(K)$, as well.
\qed
\vskip.2cm

\vskip.2cm
\noindent{\bf{Proof of Theorem~\ref{thm:ConcordanceInvariance}}}\qua
If $K$ is a slice knot (i.e.\ it bounds a
smoothly embedded disk in the four-ball), then $\tau(K)=0$, according
to Corollary~\ref{cor:GenusBounds}. The theorem follows
from this fact, together with Proposition~\ref{prop:Additivity}.
\qed
\vskip.2cm

\vskip.2cm
\noindent{\bf{Proof of Corollary~\ref{cor:SkeinInequality}}}\qua
The knot $K=K_+\#(-K_-)$ clearly bounds an immersed disk in $B^4$
with a single double-point. By resolving this double-point, we obtain
a smoothly embedded surface with genus $1$ in $W=B^4$ which bounds
$K_+\# (-K_-)$, and hence, applying
Theorems~\ref{thm:ConcordanceInvariance} and~\ref{thm:GenusBounds},
$$\tau(K_+)-\tau(K_-)\leq 1.$$ This proves one of the two
inequalities.  For the other inequality, we consider $K_-\#(-K_+)$,
and observe that if we blow up the self-intersection, we obtain an
embedded disk in $W=B^4\#\mCP$ which represents the trivial homology
class.  Thus, according to Theorems~\ref{thm:ConcordanceInvariance}
and~\ref{thm:GenusBounds}, 
$\tau(K_-)-\tau(K_+)\leq 0$.
\qed 

\section{Some small examples}
\label{sec:Examples}

In this section, we calculate $\tau$ for several ten-crossing knots,
the knots $10_{139}$, $10_{152}$, and $10_{161}$ from Rolfsen's
list~\cite{Rolfsen}. For these knots, the invariant $\tau$ gives sharp
lower bounds on the unknotting number and hence the four-ball
genus. The four-ball genera of $10_{139}$ and $10_{152}$ were first
calculated in~\cite{Kawamura} and the four-ball genus
of $10_{161}$ was first calculated 
by~\cite{Tanaka}; both results 
use gauge theory techniques, cf.~\cite{KMMilnor}
and~\cite{Rudolph}.

The present calculations rest on the techniques from~\cite{AltKnots}
and their refinements from~\cite{calcKT}. These results interpret the
generators for $\CFKa$ in terms of ``essential Kauffmann states'' for
a knot projection. We recall the definitions presently.

Let $K$ be an oriented knot, and let $G$ denote a generic projection
for $K$, with distinguished edge $\epsilon_0$. This choice of data is
called a {\em decorated knot projection} $G$. A {\em Kauffman state}
(cf.~\cite{Kauffman}) is an assignment which maps each crossing for
the knot projection $G$ one of its four adjoining quadrants, so that
no two crossings are assigned quadrants from the same region in
$S^2-G$, and no crossing point is associated to one of the two
distinguished regions containing $\epsilon_0$.

The chain complex $\CFKa$ is generated by all Kauffman states for the
knot projection.  Indeed, there is a smaller complex which can be used
when $\epsilon_0$ is chosen carefully.

To describe this, we use the notion of an {\em essential interval},
cf.~\cite{calcKT}.  Indeed, it suffices here to consider a slightly
weaker notion, which we call a {\em weakly essential interval}. A
weakly essential interval is a sequence of consecutive edges
$E=\bigcup_{i=-\ell}^m\epsilon_i$ with $\ell,m\geq 0$, so that the following
three properties hold:
\begin{itemize}
\item 
$E_-\cup\epsilon_0\cup E_+$ is an embedded arc, 
\item
as we traverse $E_+$,
all the crossings encountered have the same type (i.e.\ they are over-
or under-passes), 
\item and similarly as we traverse $E_-$, all the
crossings encountered have the same type (which might be different
from the type encountered along $E_+$).
\end{itemize}

If $v$ is a vertex in $E$ (i.e.\ a crossing for $G$), then there are
two edges in $E$ which meet $v$: of these two, one is farther from
$\epsilon_0$ (in $E$).  An {\em $E$-essential state} is a state which
associates to each vertex $v$ in $E$ one of the two regions
containing the edge through $v$ which is farther (in $E$) from $\epsilon_0$. 

\begin{figure}[ht!]
\cl{\epsfxsize 7cm\epsfbox{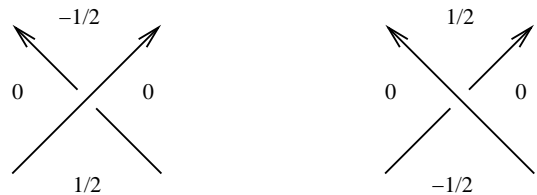}}
\caption{\label{fig:FiltLevel}
{\bf{Local filtration level contributions}}\qua We have
illustrated the local contributions for the
filtration level of a state for both kinds of
crossings.}
\end{figure}

\begin{figure}[ht!]
\cl{\epsfxsize 7cm\epsfbox{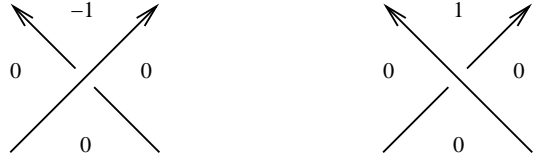}}
\caption{\label{fig:AbsGrading}
{\bf{Local grading contributions}}\qua We have
illustrated the local contribution for the absolute grading
associated for a state.}
\end{figure}

A Kauffman state $x$ is assigned a {\em filtration level} and an {\em
absolute grading}, according to the rules specified in
Figures~\ref{fig:FiltLevel} and~\ref{fig:AbsGrading} respectively. 

Finally, if we number all the edges of $G$ consecutively $\{\epsilon_i\}_{i=0}^{N-1}$,
we can associate to each state $x$ a multi-filtration-level
$$M_x\in\Hom(\{\epsilon_i\}_{i=0}^{N-1},\Z\oplus\Z)$$
by the inductive rules:
$$M_x(\epsilon_i)=
\left\{
\begin{array}{ll}
(0,0) & {\text{if $i=0$}} \\
M_x(\epsilon_{i-1})+(0,1) & {\text{if $v_i$ is over and $x(v_i)$
is to the right}}\\
M_x(\epsilon_{i-1})-(0,1) & {\text{if $v_i$ is over and $x(v_i)$
is to the left }} \\
M_x(\epsilon_{i-1})+(1,0) & {\text{if $v_i$ is an under and $x(v_i)$
is to the left }} \\
M_x(\epsilon_{i-1})-(1,0) & {\text{if $v_i$ is an under and $x(v_i)$
is to the right}} 
\end{array}\right.$$
(Here we have abbreviated the conditions.  Thus
``if $v_i$ is over and $x(v_i)$
is to the right'' means ``if $v_i$ is an overcrossing and $x(v_i)$
is to the right of $\epsilon_{i}\cup \epsilon_{i-1}$'' and the other
conditions should be expanded similarly.)
 
The following is a combination of results from~\cite{AltKnots} and~\cite{calcKT}:

\begin{theorem}
Let $G$ be a decorated knot projection for $K$, and a compatible
essential interval $E$. Then, there is a chain complex which
calculates $\CFa(S^3)$ whose generators are $E$-essential states,
whose filtration levels and absolute gradings are given as above; so,
in particular, $\HFKa(S^3,i)$ is generated by those $E$-essential
states with filtration level $i$. Moreover, the differential on
$\CFKa(S^3)$ respects the multi-filtration $M$ defined above, in the
sense that if $y$ appears with non-zero multiplicity in $\partial x$,
then for each edge $\epsilon$ not in $E$,
$M_x(\epsilon)-M_y(\epsilon)$ is a pair of non-zero integers.
\end{theorem}

\begin{proof}
In~\cite{AltKnots}, we describe a Heegaard diagram belonging to the
knot projection, for which the generators can be interpreted as
Kauffman states, with filtration level and absolute grading calculated
above. Indeed, this interpretation is established in
Theorem~\ref{AltKnots:thm:States} of~\cite{AltKnots}, while the
restriction to only essential states is described in
Proposition~\ref{KT:prop:Simplify} of~\cite{calcKT}, with a slight
modification of the original Heegaard diagram. Moreover,
compatibility with the multi-filtration is established in
Proposition~\ref{KT:prop:Domains} of~\cite{calcKT}.
\end{proof}

In the decorated knot projections we consider here, there will be a
unique maximal essential interval (through $\epsilon_0$), which we
will use as $E$.  With this understood, we will drop $E$ from our
notation.

\subsection{The knot $10_{139}$}

\begin{prop}
\label{prop:10s139}
For the knot $K=10_{139}$, we have that
$g^*(K)=u(K)=|\tau(K)|=4$.
\end{prop}

\begin{figure}[ht!]
\cl{\epsfxsize6.5cm\epsfbox{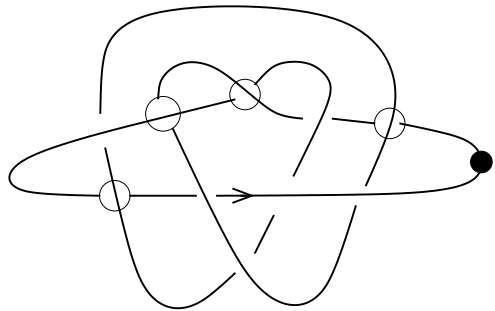}}
\caption{\label{fig:10s139}
{\bf{The knot $10_{139}$}}\qua We have circled four crossings: when these are switched,
the new knot is the unknot. The black dot indicates the distinguished edge ($\epsilon_0$)
used in the decorated knot projection.}
\end{figure}

\begin{proof}
In Figure~\ref{fig:10s139}, we have given a picture for $K=10_{139}$,
circling four crossings.  If these four crossings are changed, we
obtain a picture of the unknot. Hence, $g^*(K)\leq u(K)\leq
4$. Indeed, if we use the indicated decoration for this knot, it is
straightforward if tedious to verify that there is a unique essential
state in dimension $0$, and it has filtration level $4$. (Our
calculations were expedited by the use of Mathematica~\cite{Mathematica}.) It follows at once that
$\tau(K)=4$. The rest now follows at once from
Corollary~\ref{cor:GenusBounds}.
\end{proof}

\subsection{The knot $10_{152}$}

\begin{prop}
\label{prop:10s152}
For the knot $K=10_{152}$, we have that
$g^*(K)=u(K)=|\tau(K)|=4$.
\end{prop}

\begin{proof}
We proceed exactly as in the proof of Proposition~\ref{prop:10s139},
only with a different picture. Again, we have indicated the four
crossings which are to be changed to obtain the unknot, and we have
indicated a distinguished edge $\epsilon_0$, with respect to which
there is only one essential state in dimension zero, and it has
filtration level $-4$.  Thus $|\tau(K)|=u(K)=g^*(K)=4$.
\end{proof}

\begin{figure}[ht!]
\cl{\epsfxsize 4.5cm\epsfbox{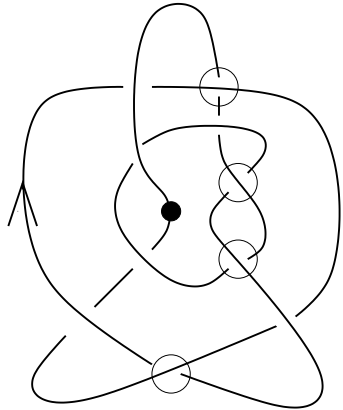}}
\caption{\label{fig:10s152}
{\bf{The knot $10_{152}$}}\qua We have indicated, once again,
the crossings needed to unknot $10_{152}$, and also 
the distinguished edge ($\epsilon_0$, with the dark circle)
used in the decorated knot projection for the proof of Proposition~\ref{prop:10s152}.}
\end{figure}

\subsection{The knot $10_{161}$}

\begin{prop}
\label{prop:10s161}
For the knot $K=10_{161}$, we have that
$g^*(K)=u(K)=|\tau(K)|=3$.
\end{prop}

\begin{figure}[ht!]
\cl{\epsfxsize 6.4cm\epsfbox{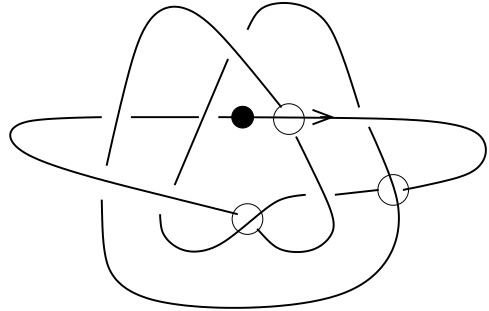}}
\caption{\label{fig:10s161}
{\bf{Local filtration level contributions}}\qua 
A decorated knot projection for $10_{161}$, showing that it has
unknotting number at least $3$.}
\end{figure}

\begin{figure}[ht!]
\cl{\epsfxsize 6.4cm\epsfbox{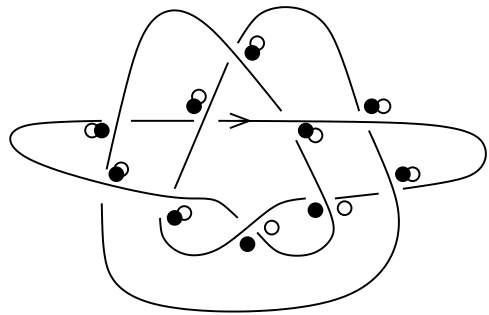}}
\caption{\label{fig:calc161}
{\bf{States $b$ and $c$}}\qua 
We have indicated the state $c$ by the dark circles, and
$b$ by the hollow ones.}
\end{figure}

\begin{figure}[ht!]
\cl{\epsfxsize 8cm\epsfbox{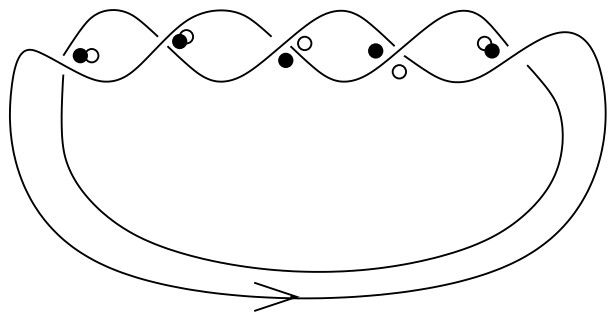}}
\caption{\label{fig:TrefCalc}
{\bf{States $y$ and $z$}}\qua 
We have indicated here the two states in absolute degrees $1$ and $2$.
The state $z$ (in dimension two) is indicated by the dark circles,
while $y$ is indicated by the hollow ones. It is easy to see that the
domain connecting $z$ to $y$ agrees with the domain connecting $c$ to $b$
in the previous picture.}
\end{figure}

\begin{proof}
First, observe that the unknotting indicated in
Figure~\ref{fig:10s161} shows that $3\leq u(K)\leq g^*(K)$. On the
other hand, for the decoration indicated in that figure, there are now
exactly two essential states with dimension zero, which we label $a$ and
$b$, where $a$ has filtration level $-3$ and $b$ has filtration
level $-2$. We claim, however, that there is an essential state $c$
in filtration level $-1$ with the property that $b$ appears once in
the expansion of $\partial c$. The states
$c$ and $b$ are illustrated in Figure~\ref{fig:calc161}.

To see that $b$ appears once in the expansion of $\partial c$, we
consider the decorated knot projection for the trefoil knot
illustrated in Figure~\ref{fig:TrefCalc}.  If we include all states
(i.e.\ we include the inessential ones), then we obtain five states,
with the extra two canceling states in filtration level $2$.  There
remain states $x$, $y$, and $z$ in filtration levels $-1$, $0$, and
$1$ respectively, and absolute degrees $0$, $1$, and $2$. It follows
at once that $\partial z=y$. On the other hand, the support of the
domain connecting $z$ to $y$ agrees with the support of the domain
connecting $c$ to $b$.

It follows now that either $c$ is null-homologous, or it is homologous
to some multiple of $b$. In either case, it follows that $b$
represents a generator of $\HFa(S^3)$, and hence $\tau(K)=-3$.
The result now follows, as usual, from Corollary~\ref{cor:GenusBounds}.
\end{proof}

\section{Generalizations to other three-manifolds}
\label{sec:Generalizations}

The constructions of the present paper can, of course, be readily
generalized to oriented links in $S^3$. Specifically, as in
Proposition~\ref{Knots:prop:LinksToKnots} of ~\cite{Knots}, an
oriented, $n$-component link $L$ in $S^3$ gives rise to an oriented
knot $\kappa(L)$ in $\#^{n-1}(S^2\times S^1)$. 
There is thus an induced filtration 
$$\Filt(\kappa(L),m)\subset \CFa(\#^{n-1}(S^2\times S^1)),$$
and hence maps
$$\iota^m_{\kappa(L)}\colon H_*(\Filt(\kappa(L),m))\longrightarrow \HFa(\#^{n-1}(S^2\times S^1))
\cong \Wedge^{*}(\Q^{n-1}).$$
Again, we define $\tau(L)$ to be the minimum $m$ for which
$\iota^m_{\kappa(L)}$ is surjective. 

In another direction, we can generalize $\tau$ to the case of
null-homologous knots in an arbitrary closed, oriented three-manifold
$Y$.  We focus presently on the case where $Y$ is an integral homology
three-sphere. There are natural maps $p\colon \HFinf(Y)\longrightarrow
\HFp(Y)$, and $q\colon \HFa(Y)\longrightarrow \HFp(Y)$. We now let
$$\tau(Y,K)=\min\{m\big| \Image(q\circ \iota^m_K)\cap \Image(p)~{\text{is non-torsion}}\}.$$

\end{document}

%% file: gtoutput.tex
\def\ifplaintex{\expandafter\ifx\csname documentclass\endcsname\relax}

\ifplaintex 
\else
\fi

\expandafter\ifx\csname beginpicture\endcsname\relax
\expandafter\ifx\csname documentclass\endcsname\relax
\input pictex \else
\input prepictex \input pictex \input postpictex \fi\fi

\def\gt{{\mathsurround=0pt\it $\cal G\mskip-2mu$eometry \&\ 
$\cal T\!\!$opology}}        

\def\gtp{{\mathsurround=0pt\it $\cal G\mskip-2mu$eometry \&\ 
$\cal T\!\!$opology $\cal P\!$ublications}}  

\def\lognumber#1{\def\thelognumber{#1}}
\def\volumenumber#1{\def\thevolumenumber{#1}}
\def\papernumber#1{\def\thepapernumber{#1}}
\def\volumeyear#1{\def\thevolumeyear{#1}}

\def\pagenumbers#1#2{\def\startpage{#1}\def\finishpage{#2}}
\def\published#1{\def\publishdate{#1}}
\def\proposed#1{\def\theproposer{#1}}
\def\seconded#1{\def\theseconders{#1}}
\def\received#1{\def\receiveddate{#1}}
\def\revised#1{\def\reviseddate{#1}}
\def\accepted#1{\def\accepteddate{#1}}

\def\coverauthors#1{\def\thecoverauthors{#1}}
\def\asciiauthors#1{\def\theasciiauthors{#1}}
\def\asciiaddress#1{\def\theasciiaddress{#1}}
\def\asciiemail#1{\def\theasciiemail{#1}}
\long\def\asciiabstract#1{\long\def\theasciiabstract{#1}}

\let\\\par\let\thelognumber\relax
\let\thevolumenumber\relax\let\thepapernumber\relax
\let\thevolumeyear\relax\let\thesamplenumber\relax\let\startpage\relax
\let\finishpage\relax\let\publishdate\relax\let\receiveddate\relax
\let\reviseddate\relax\let\accepteddate\relax\let\theasciititle\relax
\let\theasciiauthors\relax\let\theasciiaddress\relax
\let\theasciiabstract\relax
\let\theasciiemail\relax\let\theshortauthors\relax\let\theshorttitle\relax
\let\thecoverauthors\relax

\long\def\maketitlep{   

\count0=\startpage

\gt\hfill      
\beginpicture
\setcoordinatesystem units <0.33truein, 0.33truein> point at 2.2 0.9
\setplotsymbol ({$\cal G$})
\plotsymbolspacing=9truept
\circulararc 315 degrees from 0 1 center at 0 0
\setplotsymbol ({$\cal T$})
\circulararc 315 degrees from 1 -1 center at 1 0
\endpicture
\break
{\small\ifx\thesamplenumber\relax 
Volume \else Sample
\fi\thevolumenumber\ (\thevolumeyear)
\startpage--\finishpage\nl
Published: \publishdate}
\vglue 0.5truein plus 0.4fil minus 0.1truein

{\parskip=0pt\leftskip 0pt plus 1fil\def\\{\par\smallskip}{\ifplaintex\large
\else\Large\fi\bf\thetitle}\par\medskip}   

\vglue 0pt plus 0.1fil 

{\parskip=0pt\leftskip 0pt plus 1fil\def\\{\par}{\sc\theauthors}
\par\medskip}

\vglue 0pt plus 0.1fil 

{\small\parskip=0pt\let\newline\\
{\leftskip 0pt plus 1fil\def\\{\par}{\sl\theaddress}\par}
\expandafter\ifx\theemail\relax    
\relax\else\vglue 5pt plus 0.02fil minus 2pt\def\\{\stdspace{\rm 
and}\stdspace} 
\cl{Email:\stdspace\tt\theemail}\fi
\ifx\theurl\relax                  
\relax\else\vglue 5pt plus 0.02fil minus 2pt\def\\{\stdspace{\rm 
and}\stdspace}
\cl{URL:\stdspace\tt\theurl}\fi\par}

\vglue 7pt plus 0.3fil minus 3pt

{\bf Abstract}
\vglue 5pt plus 0.1fil minus 2pt

\theabstract

\vglue 7pt plus 0.3fil minus 3pt

{\bf AMS Classification numbers}\quad Primary:\quad \theprimaryclass

Secondary:\quad \thesecondaryclass

\vglue 5pt plus 0.3fil minus 2pt

{\bf Keywords}\quad \thekeywords

\vglue 10pt plus 0.5fil minus 5pt

{\small  Proposed: \theproposer\hfill Received: \receiveddate\nl
Seconded: \theseconders\hfill 
\ifx\reviseddate\relax                         
Accepted: \accepteddate                        
\else
Revised: \reviseddate                          
\fi}
\eject
}       

\let\maketitlepage\maketitlep
\let\maketitle\maketitlepage

\font\phead=cmsl9 scaled 950
\font=cmsl9 scaled 1050
\font\pnum=cmbx10 scaled 913
\font\lnum=cmbx10 
\font\pfoot=cmsl9 scaled 950
\font=cmsl9 scaled 1050
\ifplaintex
\headline{\vbox to 0pt{\vskip -4.5mm\line{\small\phead\ifnum
\count0=\startpage ISSN 1364-0380 (on line)
1465-3060 (printed) \hfill {\pnum\folio}\else\ifodd\count0\def\\{ }%
\ifx\theshorttitle\relax\thetitle\else\theshorttitle\fi\hfill{\pnum\folio}
\else\def\\{ and }{\pnum\folio}\hfill\ifx\theshortauthors\relax\theauthors
\else\theshortauthors\fi\fi\fi}\vss}}
\footline{\vbox to 0pt{\vglue 0mm\line{\small\pfoot\ifnum\count0=\startpage
\copyright\ \gtp\hfill\else
\gt, Volume \thevolumenumber\ (\thevolumeyear)\hfill\fi}\vss
}}
\else
\makeatletter
\def\@oddhead{{\small\lhead\ifnum\count0=\startpage ISSN 1364-0380 (on line)
1465-3060 (printed) \hfill {\lnum\number\count0}\else\ifodd\count0
\def\\{ }\ifx\theshorttitle\relax \thetitle \else\theshorttitle\fi\hfill
{\lnum\number\count0}\else\def\\{ and }{\lnum\number\count0}
\hfill\ifx\theshortauthors\relax 
\theauthors\else\theshortauthors\fi\fi\fi}}\def\@evenhead{\@oddhead}
\def\@oddfoot{\small\lfoot\ifnum\count0=\startpage\copyright\ \gtp\hfill\else
\gt, Volume \thevolumenumber\ (\thevolumeyear)\hfill\fi}
\def\@evenfoot{\@oddfoot}
\makeatother
\fi

\newwrite\gtoutfile
\long\gdef\makeheadfile{  
{\def\\{, }\def\s{ }
\immediate\openout\gtoutfile head.xxx
\immediate\write\gtoutfile{Proxy-for: \ifx\theasciiauthors\relax
\theauthors\else\theasciiauthors\fi\s<\ifx\theasciiemail\relax\theemail\else\theasciiemail\fi>}
\immediate\write\gtoutfile{\noexpand\\}
\immediate\write\gtoutfile{Authors: \ifx\theasciiauthors\relax
\theauthors\else\theasciiauthors\fi}
{\def\\{ }\immediate\write\gtoutfile{Title: \ifx\theasciititle\relax
\thetitle\else\theasciititle\fi}}
\immediate\write\gtoutfile{Subj-class: GT or SG or MG etc}
\immediate\write\gtoutfile{MSC-class: \theprimaryclass\ifx\thesecondaryclass\relax\else, \thesecondaryclass\fi}
\immediate\write\gtoutfile{Journal-ref: Geom. Topol. \thevolumenumber
(\thevolumeyear) \startpage-\finishpage}
\immediate\write\gtoutfile{Comments: Published by Geometry and Topology at}
\immediate\write\gtoutfile{\s\s http://www.maths.warwick.ac.uk/gt/GTVol\thevolumenumber/paper\thepapernumber.abs.html}
\immediate\write\gtoutfile{\noexpand\\}
\immediate\write\gtoutfile{}
\ifx\theasciiabstract\relax
\immediate\write\gtoutfile{\theabstract}\else
\immediate\write\gtoutfile{\theasciiabstract}\fi
\immediate\write\gtoutfile{}
\immediate\write\gtoutfile{\noexpand\\}
\immediate\write\gtoutfile{}
\immediate\closeout\gtoutfile}}  

\def\maketitlepage{\maketitlep\makeheadfile}
\let\maketitle\maketitlepage